\RequirePackage{fix-cm}

\documentclass[]{amsart}
\usepackage{graphicx}
\usepackage{amsmath,amssymb,amsthm}
\usepackage{graphicx}
\usepackage{bm}

\newtheorem{theorem}{Theorem}

\DeclareMathOperator{\sgn}{sgn}
\DeclareMathOperator{\arcsinh}{arcsinh}
\DeclareMathOperator{\arccosh}{arccosh}
\DeclareMathOperator{\arcsech}{arcsech}
\DeclareMathOperator{\am}{am}
\DeclareMathOperator{\dn}{dn}
\DeclareMathOperator{\sn}{sn}
\DeclareMathOperator{\cn}{cn}

\begin{document}
\title[Decomoposition and conformal mapping techniques]{Decomposition and conformal mapping techniques for the quadrature of nearly singular integrals}

\author[Mitchell et al.]{William Mitchell         \and
        Abbie Natkin \and
        Paige Robertson \and
        Marika Sullivan \and
        Xuechen Yu \and
        Chenxin Zhu}
\address{Macalester College\\
              1600 Grand Ave \\
              Saint Paul, MN, 55105 USA\\}
\email{wmitchel@macalester.edu}           

\date{Received: date / Accepted: date}

\maketitle

\begin{abstract}
Gauss-Legendre quadrature and the trapezoidal rule are powerful tools for numerical integration of analytic functions. For nearly singular problems, however, these standard methods become unacceptably slow. We discuss and generalize some existing methods for improving on these schemes when the location of the nearby singularity is known. We conclude with an application to some nearly singular surface integrals of viscous flow. 
\end{abstract}

\section{Introduction}
\label{sec:intro}
It is usually easy to compute definite integrals when the integrand is analytic. Gauss-Legendre and Clenshaw-Curtis quadratures generally converge rapidly for aperiodic problems, while the humble trapezoid rule is equally powerful for measuring the area under one period of a periodic integrand. 
Unhappily, this is not the case for the class of problems known as \emph{nearly singular integrals,} wherein the integrand fails to be analytic somewhere near but not on the interval of integration in the complex plane. The trapezoidal and Gauss-Legendre rules still provide geometric convergence -- that is, the logarithm of the error decreases linearly with the number of quadrature nodes -- but with an unacceptably small slope. Standard theorems relate this slope to the domain of analyticity of the integrand. In this paper we survey the existing methods of accelerating the convergence of these standard methods and we introduce several new ones, assuming that the location of the nearby singularity is known, but without using any additional information on the nature of the singularity. 

While our motivation comes from the quadrature problem in the context of solving integral equations, some very similar issues also arise in the context of interpolation from samples of a nearly singular function. This is an important component of spectral methods for differential equations when the desired solution has abrupt fronts or peaks. In fact, some of the same acceleration techniques have been independently discovered by researchers in the two communities. For example,  Johnston and Elliot's hyperbolic sine transformation \cite{johnston2005sinh} was also discovered by Tee and Trefethen \cite{tee2006rational}, while Jafari's transformation  \cite{jafari2015new} is very similar to the repeated sinh method of \cite{elliott2008iterated}. One of our aims is to assemble the results from these two communities in one place. 

We focus on two families of acceleration strategies. The first group of methods split the domain into carefully chosen subintervals and then solve subproblems on each of them. The other strategy that we consider is to change variables with a complex-analytic transformation so that the singularity lies farther away. 

An early example of a splitting method for aperiodic problems was given by Ma and Kamiya \cite{ma2002distance}, who subdivide at the real part of the singularity, assuming that this lies within the integration interval and the problem is aperiodic (in fact they subsequently use an exponential change of variables for each of the new subproblems, thereby combining both of the principal strategies considered here). 
Subsequently, Driscoll and Weideman \cite{driscoll2014optimal} gave a formula for the optimal splitting location in terms of the location of the singularity, again for the aperiodic case. Their formula is especially useful in cases where the singularity is near the endpoint of the interval in the complex plane or on the real line outside the interval, where simply using the real part of the singularity is ineffective or impossible. We develop an analogous method for periodic problems, replacing the trapezoid rule for a full period with individual Gauss-Legendre integrations on two subintervals of different sizes. 

Splitting methods are extremely simple and, as we demonstrate, can yield dramatic improvements in the convergence rate. However, we can find even better convergence rates using methods that avoid dividing the available quadrature nodes among two or more subproblems. There are a vast array of possibilities for conformal mappings that distort the neighborhood of a real interval while fixing the endpoints and remaining real-valued and monotone within the interval. The Jacobi elliptic functions are a powerful tool for designing transformations that use all of the analyticity we have assumed; we list existing methods for periodic and aperiodic problems and give a new version for the aperiodic case when the singularity lies on the real line outside the integration interval. This is a small extension of results by Trefethen, Tee and Hale \cite{tee2006adaptive,tee2006rational,hale2008new,hale2009conformal}. However, the methods employing elliptic functions are less robust than some elementary alternatives when the integrand has other challenging features like additional distant singularities or rapid growth away from the real line. We therefore list or introduce some good elementary alternatives such as the sinh transformation for aperiodic problems \cite{johnston2005sinh,tee2006rational}, our iterated sine map for periodic problems, and our quadratic transformation for aperiodic problems with real singularity.

The organization of the paper is as follows. We state and comment on the standard theorems describing the convergence of the trapezoidal and Gauss-Legendre rules in 
Sec. \ref{sec:thms}. We consider periodic problems in Sec. \ref{sec:periodic}. For aperiodic problems we first consider singularities off the real line in Sec. \ref{sec:apc} and then singularities that occur on the real line (but outside the interval of integration) in Sec. \ref{sec:apr}. We test all of the methods on a suite of integrands with various properties in Sec. \ref{sec:examples}. As an application, we then compute some nearly singular surface integrals arising in Stokes flow in Sec. \ref{sec:SLP}, followed by concluding remarks. 

\section{Discussion of the standard theorems}
\label{sec:thms}
For a periodic integrand, the convergence rate of the trapezoid rule depends on the distance from the singularity to the real line. 
\begin{theorem}[Geometric convergence of trapezoid rule \cite{trefethen2014exponentially}]
\label{thm:trap}
Suppose that a $2\pi$-periodic 
function $f$ 
is analytic
and 
satisfies  $\lvert f(z)\rvert<M$ 
within the
strip $\lvert\Im(z)\rvert <\lambda$. 
Then the difference between the integral $I = \int_{0}^{2\pi} f(x)\,dx $ and its $n$-point trapezoidal rule approximation $I_n =(2\pi/n)\sum_{j=1}^{n} f\left(2\pi j / n\right) $ satisfies
\begin{equation}
    \left\lvert I-I_n \right\rvert
    \le 
    \frac{4\pi M}{\exp(\lambda n)-1}.
\end{equation}
\end{theorem}

A similar theorem holds for Gauss-Legendre quadrature, with an ellipse replacing the infinite strip. 
\begin{theorem}[Geometric convergence of Gauss-Legendre quadrature \cite{trefethen2019approximation}] \label{thm:gl}
 Suppose that $f$ is analytic with $\lvert f(z)\rvert < M$ on the interior of the Bernstein ellipse $E_\rho$, whose foci are $\pm 1$ and whose semimajor and semiminor axis lengths sum to $\rho>1$. Let $C$ be the constant $C=64\rho^2/({15(1-\rho^{-2})})$, let $I = \int_{-1}^1 f(x)\,dx$ be the exact value of the integral, and let $I_n$ be its $n$-point Gauss-Legendre rule approximation.  Then
\begin{equation}
    \lvert I-I_n \rvert \le \frac{CM}{\rho^{2n}}.
\end{equation}
\end{theorem}
To employ this result in practice, we will often need to find the value of the ellipse parameter $\rho$ so that the ellipse passes through a given point $z\in\mathbb{C} \setminus [-1,1]$. 
We quote af Klinteberg and Barnett's useful formula \cite{af2021accurate}, 
\begin{equation}
    \rho = \rho(z) = \lvert z \pm \sqrt{z^2-1}\rvert ,\quad \textnormal{with sign chosen so that }\rho>1.
    \label{eq:rhofromz}
\end{equation}

When we apply the preceding theorems, the twin goals of maximizing $\lambda$ or $\rho$ while minimizing $M$ are in tension. If $f$ is entire or has only distant singularities and $n$ is fixed, this leads to an optimization problem for the value of $\lambda$ or $\rho$ that will produce the strongest statement from the theorem. Of course, the relative importance of $M$ decreases as $n$ grows. 
In the nearly singular case, where $\lambda\approx0$ or $\rho\approx1$, our priority is to improve the convergence rate even if this results in a large increase in $M$.  
In this paper we are assuming knowledge about the locations of the singularities of the integrand, so $\lambda$ and $\rho$ are known. However, we make no assumptions about the nature of those singularities or about the growth of $\lvert f(z)\rvert$ away from the integration interval, so $M$ is unavailable. 
Therefore we present a range of options instead of searching for an optimal strategy. 
The relatively cautious methods we consider do not make use of all of the analyticity we have assumed for $f(z)$, and are therefore less vulnerable to the danger of fast growth in $\lvert f(z) \rvert$ away from the integration interval. In contrast, the more aggressive methods make use of all of the analyticity we have assumed (these generally involve the Jacobi elliptic functions, as we will see below). The aggressive methods will boast better theoretical convergence rates, although the errors may reach machine precision before they actually decrease at the advertised rate. In contrast, the more cautious methods have (slightly) smaller convergence rates but are more likely to actually achieve these rates for challenging integrands.


\section{Methods for periodic problems}
\label{sec:periodic}
We begin with the case where the integrand is $2\pi$-periodic. We suppose that $f$ is analytic except at the points $x=2\pi k \pm Bi$ for $k\in\mathbb{Z}$, or along branch cuts extending vertically from these points to infinity. The ordinary trapezoidal rule will be an effective choice if $B$ is large, since Theorem 1 predicts errors of size $e^{-Bn}$ with $n$ function evaluations. We therefore assume that $B$ is small but nonzero, so the integral is nearly singular, and we consider several methods to use our knowledge of the domain of analyticity of $f$ in order to improve on the performance of the trapezoid rule. 
We begin with a decomposition method and then continue with several conformal mappings. Numerical examples demonstrating these techniques appear in Sec. \ref{sec:examples}. 

\subsection{Subdivision}
\label{sec:psub}
One way to integrate over a full period of $f$ is to choose an appropriately small $\delta>0$ and then split the integral into subproblems on $[-\delta, \delta]$ and $[\delta, 2\pi-\delta]$. The subproblems are not periodic, so we apply Gauss-Legendre integration for each of them.
Rescaling the subintegrals linearly to $[-1,1]$, we have 
\begin{align}
    \int_{-\pi}^\pi f(x)\,dx &= 
    \delta \int_{-1}^1 f(\delta t)\,dt + (\pi - \delta) \int_{-1}^1 f\left( \pi + (\pi-\delta)w\right)\,dw.
    \label{eq:per2glsub}
\end{align}
The two subintegrals are singular at $t=Bi/\delta $ and at $w=\frac{B i \pm \pi }{\pi-\delta}$, respectively. If we want the overall integration procedure to converge as quickly as possible and we assume $n$ is large, we should choose $\delta$ so that both subproblems have the same asymptotic convergence rate. Equivalently, both $Bi/\delta$ and $\frac{B i \pm \pi }{\pi-\delta}$ should lie on the same ellipse with foci $\pm1$ in the complex plane. Letting $m$ denote the semiminor axis length, we have the equations 
\[\frac{0^2}{1+m^2} + \frac{(B/\delta)^2}{m^2} = 1,\qquad \frac{(\pi/(\pi-\delta))^2}{1+m^2} + \frac{(B/(\pi-\delta))^2}{m^2} = 1.
\]
After some simplification we arrive at the cubic equation
\(
0 = 2\delta^3 + 2B^2 \delta - B^2 \pi.
\)
This has a unique real solution which we write in terms of hyperbolic sines as follows:
\begin{equation}
    \delta = \frac{2B}{\sqrt{3}}\sinh\left(\frac13 \arcsinh\left(\frac{3\pi\sqrt{3}}{4B}\right)\right).
    \label{eq:perdelt}
\end{equation}
The error of the trapezoid rule on the original integral decays like $e^{-Bn}$ for large $n$, where $n$ is the number of quadrature nodes. 
With subdivision as in \eqref{eq:per2glsub} and Gauss-Legendre quadrature with $n/2$ nodes on each subinterval,\footnote{We experimented with unequal distributions of nodes between the two subintervals but did not see more than modest improvements. For odd $n$, one should assign the extra node to the larger subinterval since errors there are multiplied by $(\pi-\delta)$ instead of $\delta$ in \eqref{eq:per2glsub}. } the total error will decay like $\rho^{-2(n/2)} = \rho^{-n}$, where $\rho = (B + \sqrt{B^2 + \delta^2}) / \delta$. Therefore, the splitting procedure is advantageous when $\log\rho > B$.  We found numerically that this holds for $B<0.95$, although in practice, for finite values of $n$, it may be advisable to use the trapezoid rule for slightly smaller values of $B$.

\subsection{Conformal maps for periodic problems}
\label{sec:permaps}
A more powerful method for periodic problems is to compute the integral using the $n$-point trapezoidal rule following a transformation $x(t)$ that preserves the interval $[-\pi,\pi]$: 
\begin{equation}
   \int_{-\pi}^\pi f(x)\,dx = \int_{-\pi}^\pi f(x(t))x'(t)\,dt \approx \frac{2\pi}{n}\sum_{j=1}^n f\left(x\left(t_j\right)\right)x'\left(t_j\right)
\end{equation}
where $t_j=-\pi + {2\pi j}/{n}$. 
In view of Theorem 1, the error will decay like $e^{-\lambda n}$ as long as the new integrand $f(x(t))x'(t)$ is periodic and analytic on the strip $S_\lambda = \{t\in\mathbb{C}:\lvert \Im (t) \rvert<\lambda\}$.  We therefore seek a transformation $x(t)$ which carries a wide strip surrounding the real line in the complex $t$-plane into the domain of analyticity of $f$ in the complex $x$-plane, avoiding the branch cuts; another consideration is that the derivative $x'(t)$ should not have any singularities for $t\in S_\lambda$.  

Two recent works have presented transformations with the general properties we seek: Tee constructed one using the Jacobi elliptic functions \cite{tee2006adaptive}, while Berrut and Elefante gave an elementary formula based on a M\"obius transformation \cite{berrut2020periodic}. After discussing these two methods, we propose a third, which we call the iterated sine map. We illustrate the conformal maps in Figure \ref{fig:cmp} with $B=0.3$ along with the resulting predicted convergence rates. We then give numerical examples comparing these three mappings, as well as the decomposition method from Section \ref{sec:psub} and the ordinary trapezoid rule, in Fig. \ref{fig:many_periodic_examples}.

\subsubsection{The Jacobi amplitude map}
Tee and Trefethen designed a conformal map that carries a strip onto the doubly slit region of analyticity of $f$ \cite{tee2006adaptive}.  Here we present a new derivation of the same formula using differential equations rather than complex analysis. Our approach is based on the intuition that the factor $x'(t)$ should be small when $x$ is close to the singularity. To respect periodicity, we  wrap the real line onto the unit circle and imagine the singularity as a nearby point $(1,0,B)$. We then assume that $x'(t)$ is proportional to the distance from $(\cos(x),\sin(x),0)$ to the singularity, leading to the boundary value problem
\begin{equation}
\frac{dx}{dt} = k \sqrt{(1-\cos( x))^2 + \sin^2( x) + B^2},\quad x(-\pi)=-\pi,\;\; x(\pi)=\pi    
\label{eq:bvp2}
\end{equation}
where the value of $k$ is determined as part of the solution. The exact solution can be obtained from the Jacobi amplitude function $\am(t,m)$, using the relation 
\begin{equation}
    \frac{d}{dt}\am(t,m) = \sqrt{1 - m \sin^2\big(\am(t,m)\big)} = \dn(t,m).
\end{equation}
The solution of the BVP \eqref{eq:bvp2} is explicitly 
\begin{equation}
x(t) = -\pi + 2 \am\left(\frac{\pi + t}{\pi}K\left(\frac{4}{4+B^2}\right),\frac{4}{4+B^2}\right)
\label{eq:JAM}
\end{equation}
where the real quarter-period $K(m)$ is defined for $m<1$ by 
\begin{equation}
    K(m) = \int_0^{\pi/2} \frac{d\theta}{\sqrt{1-m\sin^2\theta}}.
\end{equation}
In practice, we suggest computing the derivative $dx/dt$ using the Jacobi elliptic function $\dn(t,m)$ rather than the differential equation \eqref{eq:bvp2}: 
\begin{equation}
    x'(t) = \frac{2}{\pi}K\left(\frac{4}{4+p_y^2}\right) \dn\left(K\left(\frac{4}{4+p_y^2}\right)(t-1),\frac{4}{4+p_y^2}\right).
\end{equation}
This leads to an improved convergence rate of 
\begin{equation}
    \lambda = \pi\frac{K\left(B^2/(4+B^2)\right)}{K\left(4/(4+B^2)\right)}.
\end{equation}
The transformation \eqref{eq:JAM} is the sum of the identity map and a  $2\pi$-periodic function, and it carries the rectangle $\{x+iy: \lvert x\rvert \le\pi, \lvert y \rvert < \lambda\}$ onto the doubly slit region $\{x+iy:\lvert x \rvert\le\pi\text{ and }\lvert y \rvert<p_y \textnormal{ if }x=0\}$, thereby using all of the analyticity that we have assumed for $f$. This is the most aggressive option for the periodic problem. 

\subsubsection{The boundary correspondence map}
Berrut and Elefante recently suggested an alternative mapping that avoids the use of elliptic functions \cite{berrut2020periodic}. Adapting their method slightly, we get the formulas 
\begin{align}
    a &= \exp(B) - \sqrt{-1+\exp(2B)}\\
    x(t) &= -i \log\left(\frac{\exp(it)+a}{1+a\exp(it)}\right)\\
    x'(t)&= 
    \frac{1-a^2}{a^2+2 a \cos (t)+1}\textnormal{ for } t\in\mathbb{R}\\
    \lambda_{\max} &= -\log(a).
\end{align}
This method is much better than the ordinary trapezoid rule, but it does not perform nearly as well as the Jacobi amplitude map, as one can see in Fig. \ref{fig:many_periodic_examples}. The conformal image of the strip (top right in Fig. \ref{fig:cmp}) is unbounded and does not closely approach the sides of the vertical branch cuts. 
\subsubsection{The iterated sine map}
\begin{figure}
    \centering
    \includegraphics[width=\linewidth]{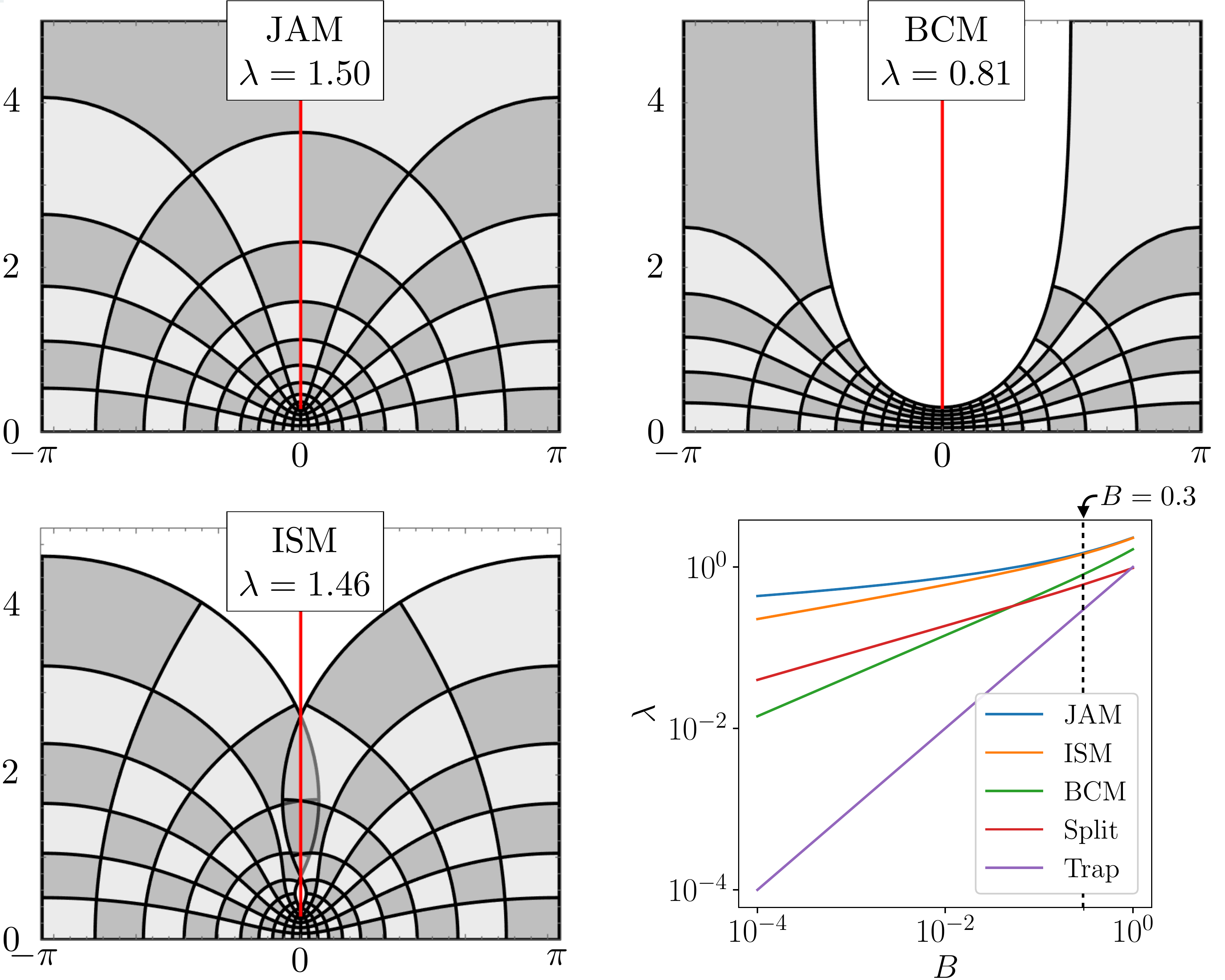}
    \caption{Comparison of three conformal mappings for periodic problems when the integrand has a branch cut extending vertically from $0+0.3i$ (red lines).  The three grids are the conformal images of three rectangles $\{x+iy: \lvert x \rvert\le\pi, 0\le y\le \lambda\}$. The Jacobi amplitude mapping carries a wide strip $(\lambda = 1.5)$ onto the full, unbounded domain of analyticity we have assumed for the integrand. The boundary correspondance map carries a thinner strip ($\lambda = 0.81$) into a subset of the region of analyticity, avoiding the sides of the branch cut but extending to infinity for $\lvert x \rvert>\pi/2$.  The iterated sine map carries a wide strip $(\lambda = 1.46$) onto a bounded region which slightly overlaps the branch cut, an issue that disappears when the height $B$ of the singularity is less than $0.03$.  At bottom right, we compare the improved convergence rates for these conformal maps as well as the splitting method and the ordinary trapezoid rule for $B\in[10^{-4},1]$.
    }
    \label{fig:cmp}
\end{figure}

We searched for another conformal mapping with the goal of reproducing the good convergence rate of the Jacobi amplitude map with a simpler transformation.
An early candidate was the mapping $\phi(t) = t - a \sin(t)$ for $a\in[0,1)$. 
This function increases monotonically, it carries $[-\pi,\pi]$ to itself, and it has a small derivative when $0=t=x(t)$ if $a\approx 1$. 
These properties also hold for the compositions $\phi\circ \phi$ and  $\phi\circ\phi\circ\phi$. 
The mapping $x(t) = \phi(t)$ is not competitive with the Jacobi amplitude transformation, but with a suitable choice of parameter $a$, the mapping  $x(t)=\phi(\phi(t))$ is. We call this the \emph{iterated sine map}: 
\begin{align}
    x(t) &= t - a\sin(t) -  a\sin(t-a\sin(t))\\
    x'(t) &= (1-a\cos t)\left(1-a\cos(t-a\sin t)\right).
\end{align}
The optimal value of the parameter $a$ depends on the location of the singularity $0\pm Bi$. 
To find it, we study the imaginary part of  $x(0+i\lambda)$ as a function of $\lambda$. This increases from $0$ when $\lambda =0$ to a local maximum when $\lambda=\arccosh(1/a)$. Setting the value of this local maximum equal to $B$, we have the equation
\begin{equation}
    -\sqrt{1-a^2}+a \sinh \left(\sqrt{1-a^2}-\arcsech(a)\right)+\arcsech(a)=B.
\end{equation}
This is difficult or impossible to solve analytically for $a$. In searching for a good initial guess for Newton's iteration, we found that the optimal value of $a$ was very close to $1 + B/5 - B^{2/5}$, especially for $B \approx 0$. The approximation is good enough that we forego Newton iteration entirely and use $a = 1+B/5-B^{2/5}$, with the caveat that for $B>1.5$, one should abandon the iterated sine map and use the ordinary trapezoid rule without any change of variable. This mapping is illustrated at bottom left in Fig. \ref{fig:cmp}; it is a `cautious' method because it carries the strip into a relatively small bounded region, but its convergence rate is nearly as good as the Jacobi amplitude map. For $B>0.03$ the image of the strip overlaps the branch cut slightly, possibly impairing the convergence rate, but this problem disappears for $B<0.03$. 

\section{Methods for aperiodic problems with singularity off the real line}
\label{sec:apc}
We now consider the problem of integrating an aperiodic function $f(x)$ on $[-1,1]$, assuming that $f$ is analytic except for singularities at $A\pm Bi$ with $B>0$, or along branch cuts extending vertically from these points to infinity. Our goal is to use this information about the integrand to improve on standard Gauss-Legendre quadrature. As before, we describe a decomposition method and also several methods based on conformal maps. Here we do not propose any novel mappings, but we give some new results on the convergence properties of those previously described \cite{tee2006rational,johnston2005sinh,elliott2008iterated,jafari2015new}. 
\subsection{Decomposition}
A simple strategy for integration on $[-1,1]$ is to carry out separate Gauss-Legendre integrations on $[-1,\delta]$ and $[\delta,1]$, where $\delta$ depends on $A\pm Bi$. Rescaling both subproblems back to $[-1,1]$, we obtain 
\begin{equation}
\footnotesize{   \int_{-1}^1 f(x)\,dx 
    = \frac{\delta+1}{2}\int_{-1}^1 f\left(\frac{\delta-1}{2} + \frac{\delta+1}{2}t \right)\,dt  + \frac{1-\delta}{2}\int_{-1}^1 f\left(\frac{1+\delta}{2} + \frac{1-\delta}{2}w \right)\,dw.}
    \label{eq:glchop}
\end{equation}
Now the integrands on the right side of \eqref{eq:glchop} have singularities at the points $t^* = (2A-\delta+1+ 2Bi)/(1+\delta)$ and $w^* = ({2A-\delta-1+2Bi})/({1-\delta})$, respectively. To optimize the convergence rate, we must choose $\delta$ so that the Bernstein ellipses passing through $t^*$ and $w^*$ coincide. In particular, their semimajor axis lengths must be equal. Letting $a$ denote this common length, we have the equations
\begin{equation}
    \frac{(2A-\delta+1)^2}{a^2(\delta+1)^2} + \frac{(2B)^2}{(a^2-1)(\delta + 1)^2} = 1 = 
        \frac{(2A-\delta-1)^2}{a^2(\delta+1)^2} + \frac{(2B)^2}{(a^2-1)(1-\delta )^2} .
\end{equation}
By eliminating $a$ we arrive at the quartic polynomial equation
\begin{equation}
    2\delta (2A - \delta+1)^2 (A-\delta) + 4\delta B^2 (2A-\delta) = 2(\delta+1)^2(2A-\delta)(A-\delta).
\end{equation}
The real solution with $\lvert \delta \rvert <1$ is given by 
\begin{equation}
\label{eq:apdelta}
\delta=\sgn(A)\left(
\lvert A \rvert  - \frac{\sqrt{2}}{2}\sqrt{A^2-1 - B^2 + \sqrt{-4A^2 + \left(1+A^2+B^2\right)^2}}\right).
\end{equation}
This formula agrees with Equation 2.7 in \cite{driscoll2014optimal}, and gives a simple method for evaluating $\int_{-1}^1 f(x)\,dx$, given the location of the nearest singularity $x^* = A+Bi$: define $\delta$ using \eqref{eq:apdelta}, then evaluate the right-hand side of \eqref{eq:glchop} using $n/2$-point Gauss-Legendre quadrature on each subintegral.  
Letting $\rho_{x^*}$ and $\rho_{t^*}=\rho_{w^*}$ be the Bernstein ellipse parameters for the original and decomposed problems, we see that the error decay rates of Theorem \ref{thm:gl} are $\rho_{x^*}^{-2n}$ and $\rho_{t^*}^n$, respectively. Therefore the splitting procedure will be worthwhile if $\sqrt{\rho_{t^*}} > \rho_{x^*}$. 

\subsection{Conformal maps} 
\begin{figure}
    \centering
    \includegraphics[width=\linewidth]{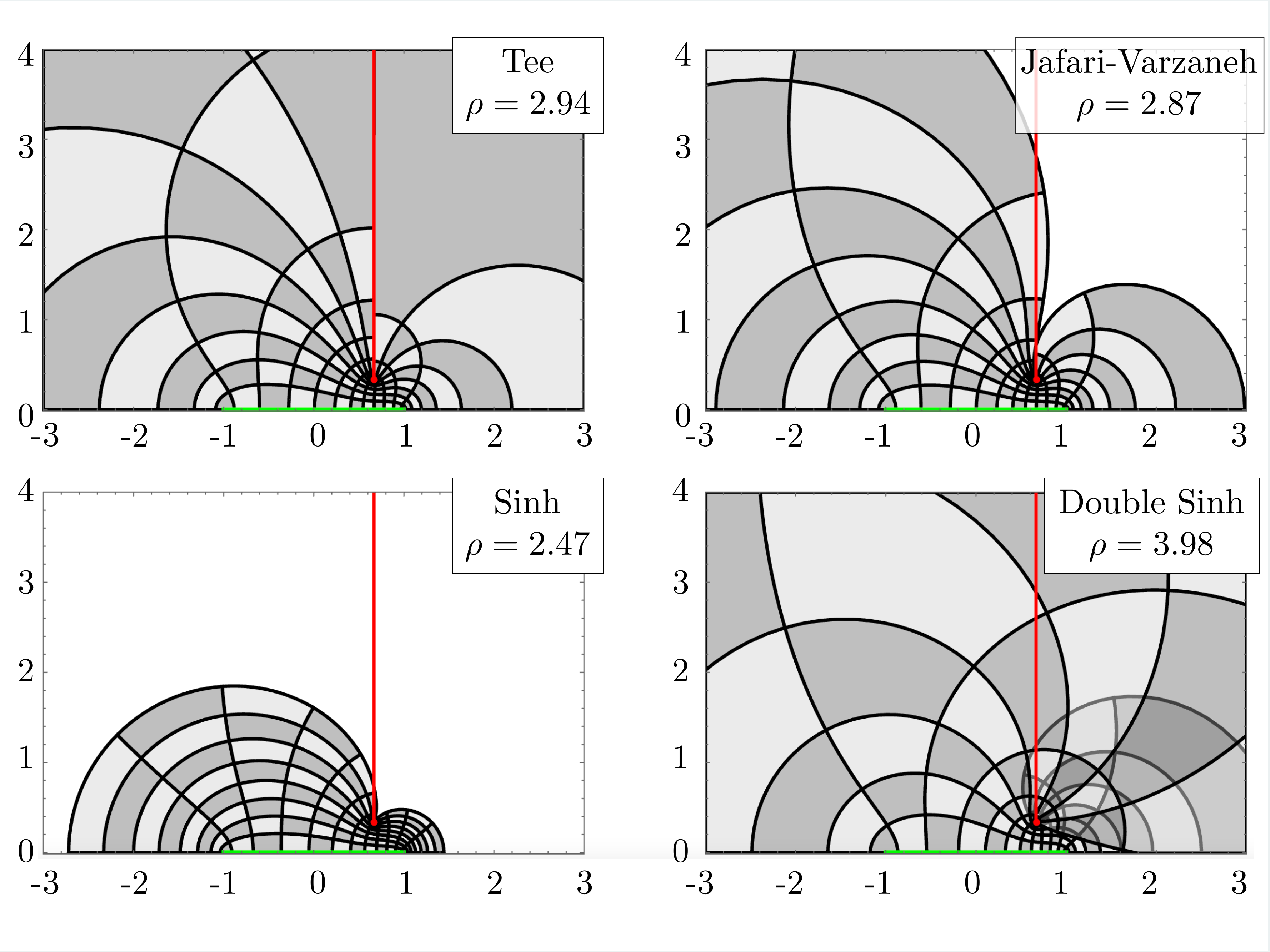}
    \caption{Comparison of four conformal mappings designed to avoid a singularity at $2/3 + i/3$. We plot the image of four Bernstein ellipses with different values of $\rho$. The mapping introduced by Tee, in terms of the Jacobi elliptic sine, permits a large value of $\rho$ and uses all of the assumed analyticity. The $\rho$-value for the map of Jafari-Varzaneh and Hosseini is nearly as large, and the image is bounded. The hyperbolic sine mapping has the smallest $\rho$-value of the four, and it uses far less of the assumed domain of analyticity.  
    The iterated or double sinh transformation achieves a large value of $\rho$, but the image of the ellipse extends far from the original interval and, while avoiding the singularity itself, wraps around it and significantly overlaps any branch cut extending from the singularity. 
    }
    \label{fig:conformalC}
\end{figure}

We now turn to strategies based on conformal mapping. Writing $t_j$ and $w_j$ for the nodes and weights of Gauss-Legendre quadrature on $[-1,1]$, we have
\begin{equation}
    \int_{-1}^1 f(x)\,dx = \int_{-1}^1 f(x(t))x'(t)\,dt \approx \sum_{j=1}^n f(x(t_j))x'(t_j)w_j
\label{eq:GLtransform11}
\end{equation}
The new integrand $f(x(t))x'(t)$ should be analytic on a large Bernstein ellipse in the complex $t$-plane. In particular, the image of the ellipse should lie within the domain of analyticity of $f$ in the complex $x$-plane, and the derivative $x'(t)$ should not itself be singular, at least for $t$ within the Bernstein ellipse. We also require that $x(t)\in[-1,1]$ for $t\in[-1,1]$, and $x(\pm1)=\pm 1$.  
See Fig. \ref{fig:conformalC} for illustrations of the three transformations when $A\pm Bi = 2/3 \pm i/3$ as well as a plot of the convergence rates for other values of $B$, and Fig. \ref{fig:many_np_c_examples} for numerical results. 
\subsubsection{The hyperbolic sine map}
We give a new derivation of the sinh transformation \cite{johnston2005sinh,tee2006adaptive}, using a differential equation. Motivated by the desire to keep the derivative $x'(t)$ small when $x$ is near the singularity $A+Bi$, we consider the BVP
\begin{equation}
    x'(t) = k\sqrt{(x-A)^2 + B^2},\quad x(-1)=-1,\;\;x(1)=1
\end{equation}
where the proportionality constant $k$ is to be found as part of the solution. 
This can be solved by hand, yielding the transformation 
\begin{align}
    x(t) &=  A + B \sinh\left( \frac{1-t}{2}\arcsinh\left(\frac{-1-A}{B}\right) + \frac{1+t}{2}\arcsinh\left(\frac{1-A}{B}\right) \right).
\end{align}
The derivative $x'(t)$ is an entire function, so the convergence rate will be limited by the analyticity of $f(x(t))$. 
The value of $t$ for which $x(t)=A+Bi$ is
\begin{equation}
    t^* = 1 + \frac{i \pi - 2 \arcsinh((1-A)/B)}{\arcsinh((1-A)/B + \arcsinh((1+A)/B}
\end{equation}
and we can use this in \eqref{eq:rhofromz} to obtain the improved value of the ellipse parameter $\rho$. 

\subsubsection{The elliptic sine map}
By conformally mapping a Bernstein ellipse onto the doubly slit plane, Tee and Hale \cite{hale2009conformal} obtained the transformation $x(t)$ given by the equations
\begin{align}
c &= \frac{\sgn(A)}{\sqrt{2}} \sqrt{A^2+B^2+1 - \sqrt{(A^2+B^2+1)^2-4A^2}}\\
m &= \frac{\left(-B + \sqrt{B^2 + 1-c^2}\right)^4}{(1-c^2)^2}\\
h(t) &= m^{1/4}  \sn\left(\frac{2K(m)}{\pi}\arcsin(t), m\right)\\
x(t)&= \frac{c}{m^{1/4}} - \frac{1-m^{1/2}}{2m^{1/4}} \left(\frac{1-c}{ h(t) - 1} + \frac{1+c}{h(t)+1}\right).
\end{align}
Note that we have condensed their notation. The intermediate conformal map $h(t)$, carrying an ellipse to the unit circle, is due to Schwarz \cite{szego1950conformal,schwarz1869}. 
The enlarged Bernstein ellipse has  parameter given by \cite{hale2009conformal}
\begin{equation}
    \rho = \exp\left(\frac{\pi K(1-m)}{4 K(m)}\right).
\end{equation}
For convenience in computing the  derivative $x'(t)$, we note that 
\begin{equation}
    h'(t) = \frac{2m^{1/4}K(m)}{\pi\sqrt{1-t^2}} \cn\left(\frac{2K(m)}{\pi}\arcsin(t), m\right)\dn\left(\frac{2K(m)}{\pi}\arcsin(t), m\right).
\end{equation}
This aggressive strategy takes full advantage of the assumed analyticity of the integrand and has a large convergence rate.  
\subsubsection{Jafari-Varzaneh's mapping}
Jafari-Varzaneh and Hosseini used the composition of the hyperbolic sine mapping and another similar function to obtain a new mapping \cite{jafari2015new}. With \(\alpha = \frac12 \arcsinh\left(\frac{1-A}{B}\right) + \frac12 \arcsinh\left(\frac{1+A}{B}\right)\), their formulas can be condensed to 
\begin{align}
    x(t) &= A + B \sinh\left(\arcsinh\left(\frac{1-A}{B}\right)-\alpha
    + \frac{2\alpha L + \pi}{2} \tan\left(t\, \arctan\left(\frac{2\alpha}{2\alpha L + \pi}\right)\right)\right)
\end{align}
where $L$ is a tunable parameter between $0.2$ and $0.9$. In light of their statement that ``experiments show that different values of this parameter approximately give the same results,'' we choose to take $L=0.5$ in all cases. 
We solved for the value of  $t^*$ satisfying $x(t^*) = A+Bi$ and obtained 
\begin{equation}
    t^* = \arctan\left(\frac{2\alpha - 2\arcsinh\left(\frac{1-A}{B}\right)+i\pi)}{2\alpha L+\pi}\right)/\arctan\left(\frac{2\alpha}{2\alpha L+\pi}\right),
    \label{eq:sinhmapsing}
\end{equation}
which can be used in \eqref{eq:rhofromz} to find a prediction of the new convergence rate. This strategy is more cautious than Tee's elliptic sine mapping but more aggressive than the hyperbolic sine tranformation. 
\subsubsection{The iterated sinh map}
For some problems, Elliot and Johnston suggested applying the sinh mapping twice in succession \cite{elliott2008iterated}. To construct the mapping, we can use \eqref{eq:sinhmapsing} to obtain the singularity location after one sinh transformation and then apply another sinh transformation accordingly. For convenience, we simplify to obtain the formula
\(x(t) = A+B \sinh \left(\frac{\pi}{2}  \sinh  \left(\ell(t)\right)\right)\), where $\ell(t)$ is the linear function
\begin{equation}\ell(t) = \frac{t+1}{2} \arcsinh\left(\frac{2}{\pi} \arcsinh\left(\frac{1-A}{B}\right)\right)+\frac{t-1}{2} \arcsinh\left(\frac{2}{\pi} \arcsinh \left(\frac{A+1}{B}\right)\right).
\end{equation}
The new singularity lies at 
\begin{equation}
t^* = 1+\frac{i \pi-2 \arcsinh\left(\frac{2}{\pi} \arcsinh((1-A)/B)\right)}{\arcsinh\left(\frac{2}{\pi} \arcsinh((1-A)/B)\right)
+\arcsinh\left(\frac{2}{\pi} \arcsinh((1+A)/B)\right)}
\label{eq:sinhsinhtstar}
\end{equation}
instead of $x^*=A\pm Bi$. 
This map tends to wrap the original Bernstein ellipse onto a relatively large region surrounding the singularity (overlapping the branch cut if there is one). Elliot and Johnston observed that the iteration gives better results for what they term \emph{nearly strongly singular} problems, meaning integrals that become divergent when the singularity reaches the real line. In contrast, they present experiments showing that a single sinh transformation is better for \emph{nearly weakly singular} integrals (which become convergent improper integrals when the singularity reaches the integration interval). This distinction is outside our scope because it relies on information about the nature of the singularity in addition to its location. 

\section{Methods for aperiodic problems with singularity on the real line}
\label{sec:apr}
We now consider the case where the integrand $f(x)$ is analytic except for an isolated singularity at some real $A$ with $\lvert A \rvert >1$, or possibly with a horizontal branch cut from $A$ to infinity. This situation is less common than the fully complex case, both for quadrature problems and for rational barycentric interpolation and Chebyshev spectral methods, and accordingly has received less attention. We describe several strategies for using the analyticity of $f$ in order to improve on Gauss-Legendre integration with $n$ nodes; to our knowledge these are all new. 

\subsection{Subdivision}
\label{sec:arsub}
One simple way to take advantage of information on the location of the singularity is to carry out separate Gauss-Legendre integrations, each using $n/2$ nodes, on $[-1,\delta]$ and $[\delta,1]$, where $\delta$ is given by 
\eqref{eq:apdelta}. With $B=0$, that formula simplifies to 
\begin{equation}
    \delta = \sgn(A) 
    \left(
    \lvert A \rvert  - \sqrt{A^2-1}
    \right).
\end{equation}
As above, the improved convergence rate comes from putting the rescaled singularity $(2A-\delta+1)/(1+\delta)$ into \eqref{eq:rhofromz}.


\subsection{Conformal maps}
\begin{figure}
    \centering
    \includegraphics[width=\linewidth]{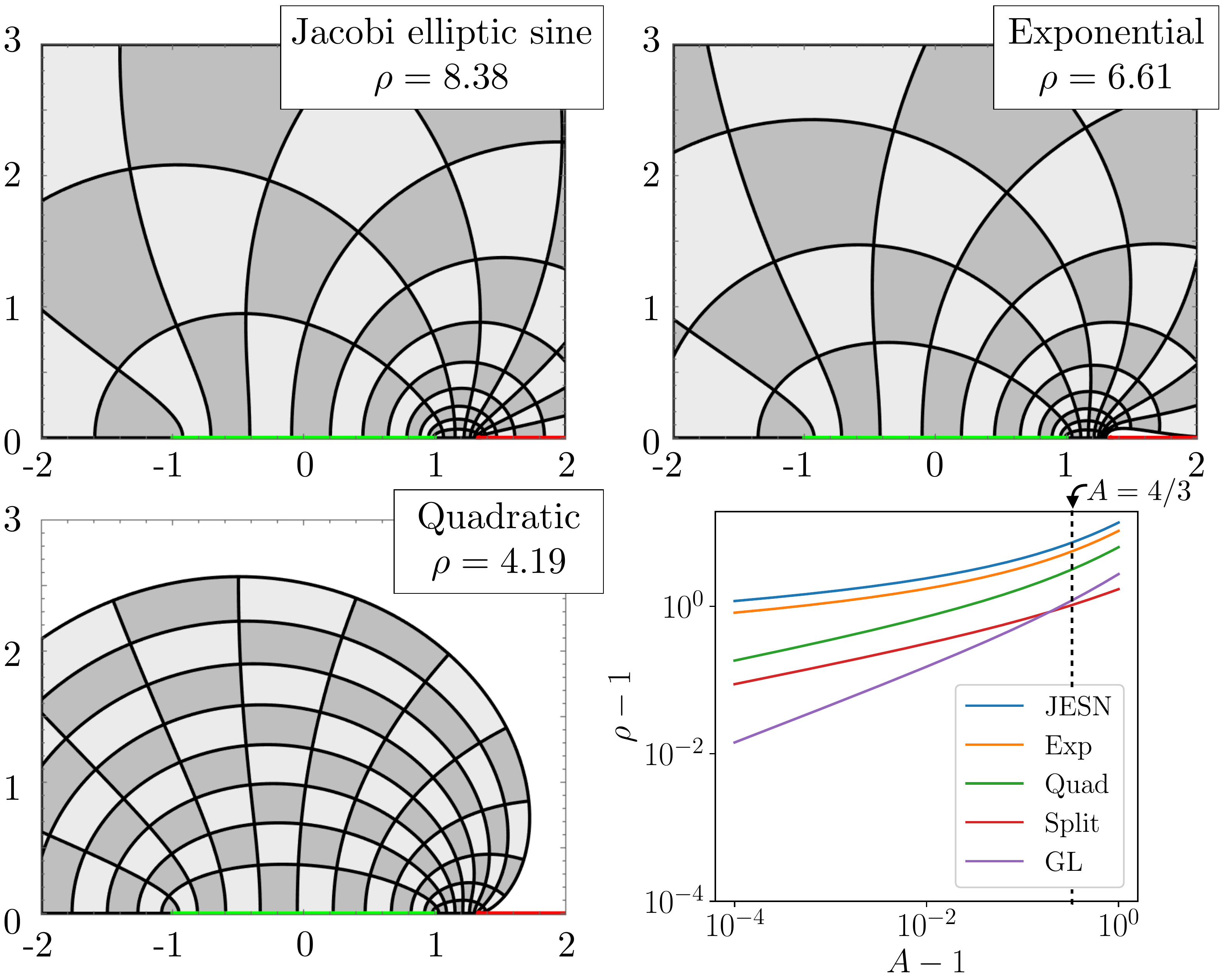}
    \caption{Three conformal mappings that avoid a singularity on the real line, at $4/3.$ The Jacobi elliptic sine mapping (top left) carries a large ellipse ($\rho = 8.38$) onto $\mathbb{C}\backslash[4/3,\infty)$. The exponential transformation (top right) carries a smaller ellipse ($\rho = 6.61)$ onto a large but bounded region. The quadratic transformation (bottom left) carries a smaller ellipse ($\rho=4.19$) into a much smaller bounded region. All three methods improve substantially against ordinary Gauss-Legendre integration, which has $\rho = 2.21$ for $A=4/3$. At lower right we plot the relationship between $\rho-1$ and $A-1$ for these methods as well as the decomposition method of section \ref{sec:arsub}. }
    \label{fig:complextransformsReal}
\end{figure}

We again follow \eqref{eq:GLtransform11}, applying Gauss-Legendre integration following a change of variable. The mapping $x(t)$ must preserve the interval $[-1,1]$ and fix its endpoints, and ideally should carry a large Bernstein ellipse in the $t$-plane into the domain of analyticity of $f$. We introduce three possibilities and illustrate them in Fig. \ref{fig:complextransformsReal}; numerical examples appear in Fig. \ref{fig:many_np_r_examples}. 

\subsubsection{Quadratic transformation}
There is a unique second-degree polynomial $x(t)$ which satisfies the boundary conditions $x(-1)=-1$ and $x(1)=1$, satisfies $x'(t)=0$ when $x=A$, and is strictly increasing for $-1\le t\le1$. 
It is given by 
\begin{equation}
    x(t) = \frac{-1}{2} \sgn(A) \left(\lvert A \rvert -\sqrt{A^2-1}\right)(t^2-1) + t.
\end{equation}
The unique value of $t$ satisfying $x(t)=A$ is the semimajor axis length of the enlarged ellipse; this leads to the improved ellipse parameter 
\begin{equation}
    \rho = \lvert A \rvert  + \sqrt{A^2-1} + \sqrt{(\lvert A \rvert  + \sqrt{A^2-1})^2-1}.
\end{equation}
This relatively cautious strategy gave the best results for most of our test problems. 
\subsubsection{Exponential transformation}
By solving the BVP 
\begin{equation}
    \frac{dx}{dt} = k\lvert x-A\rvert ,\quad x(-1)=-1,\;x(1)=1,
\end{equation}
where the value of $k$ is determined as part of the problem, we obtain an exponential change of variable:
\begin{equation}
    x(t) = A + (1-A)\exp\left(\frac{1-t}{2} \log\left(\frac{A+1}{A-1}\right)\right)
\end{equation}
There is no value of $t$ satisfying $x(t)=A$, so the composition $f(x(t))$ is entire if $f$ has an isolated singularity. On the other hand, if $f$ has a branch cut from $A$ to infinity, then the composition is analytic only on the Bernstein ellipse with parameter $\rho = s + \sqrt{1+s^2}$ with 
$s = 2\pi/\lvert\log( (A+1)/(A-1))\rvert.$  
\subsubsection{Elliptic sine map}
We have already seen the Schwarz mapping carrying the interior of an ellipse onto the unit disk. We can compose this map with the transformation 
$z\mapsto  - (1-z)^2/(1+z)^2$ to obtain a mapping from the ellipse to the slit plane $\mathbb{C}\backslash[0,\infty)$.
Then, by scaling and translation, we can arrange for $x(\pm 1) = \pm 1$. The transformation and its derivative are given by
\begin{align}
c_1 &= 17-80A^2+64A^4\\
c_2 &=  (4A^2-3)\lvert A \rvert \sqrt{A^2-1}\\
m &= c_1+16c_2 - 4 \sqrt{2}\sqrt{(A^2-1)(8+c_1(4A^2-1)) + c_1c_2}\\
h(t) &= m^{1/4}  \sn\left(\frac{2K(m)}{\pi}\arcsin(t), m\right)\\
x(t)&= \frac{h(t)+2\sqrt{m}\left(1+h(t)+h(t)^2\right) + m\, h(t)}{(1+h(t))^2 \left(1+\sqrt{m}\right)m^{1/4}}\\
x'(t) &= -h'(t) \frac{(h(t)-1)\left(1-\sqrt{m}\right)^2}{(1+h(t))^3\left(1+\sqrt{m}\right)m^{1/4}}
\end{align}
while the new convergence rate is $\rho = \exp(0.25\pi K(1-m)/K(m))$. This aggressive strategy was more effective than the splitting method but less effective than the other conformal maps in our tests. 

\section{Numerical examples}
\label{sec:examples}
We now assess the performance of all of the strategies described above on a collection of nearly singular integrals of varying difficulty. We begin with a careful discussion of the results for periodic problems and then treat the  aperiodic problems more briefly, since the themes are similar.  

\subsection{Periodic examples} 
\begin{figure}
    \centering
    \includegraphics[width=\linewidth]{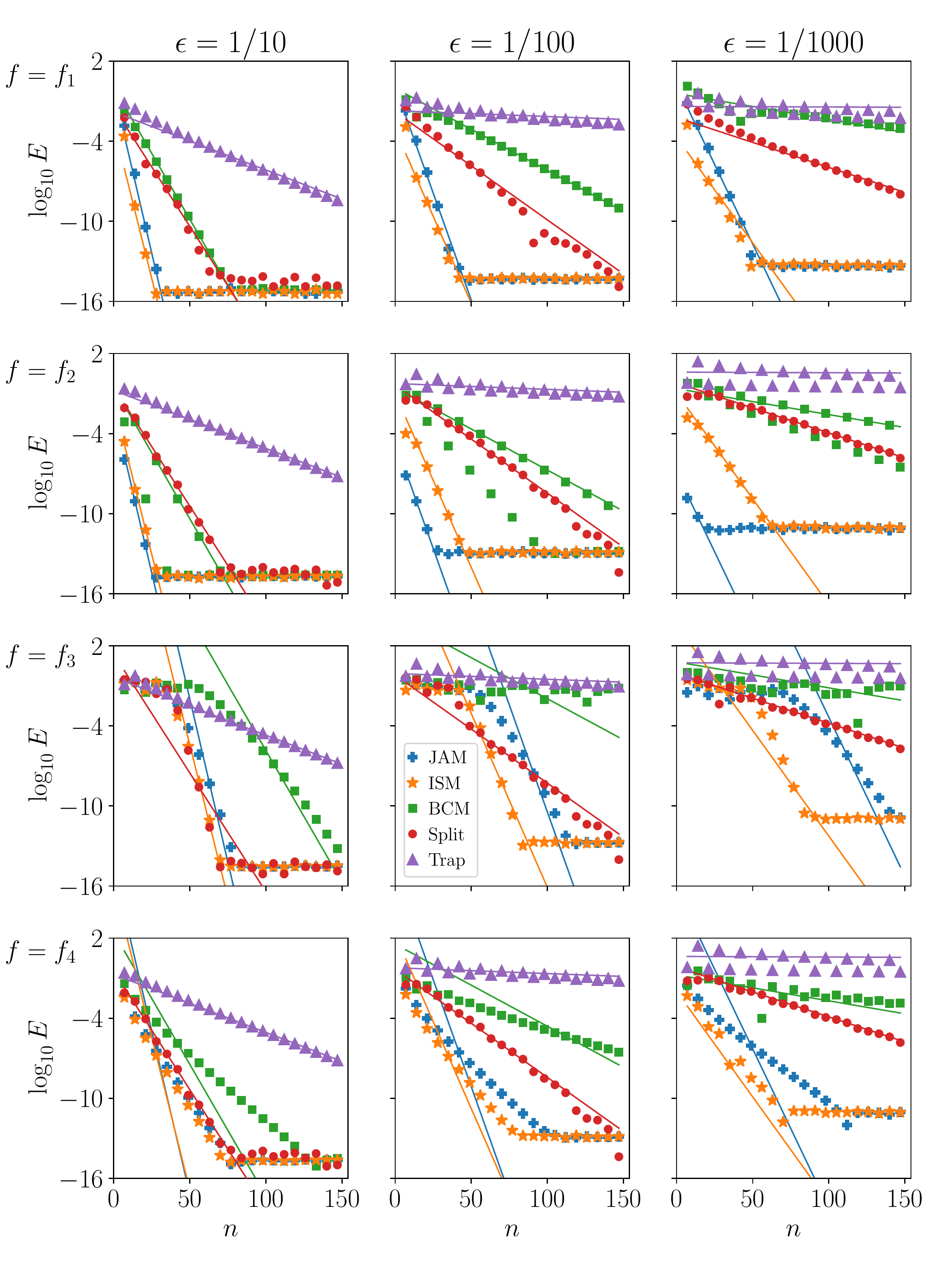}
    \caption{We test five strategies for integration of periodic, nearly singular integrals on twelve problems of varying difficulty. The relative numerical errors (glyphs) generally decay with the predicted slopes (lines), with exceptions in the last two rows due to particular features of the integrands as discussed in the text. The predicted slopes depend only on the locations of the singularities, $\pm \epsilon i$, and are identical within columns of the figure. The iterated sine mapping (stars) gives the best results, except in the second row where the Jacobi amplitude transformation (pluses) is better. }
    \label{fig:many_periodic_examples}
\end{figure}

For periodic problems, we consider integrals of the form  $\int_{-\pi}^\pi f_i(x;\epsilon)\,dx$ for $\epsilon\in\{10^{-1}$, $10^{-2}$, $10^{-3}\}$ where $f_i$ is one of the following:
\begin{align}
    f_1(x;\epsilon) &= \log(\cosh(\epsilon)-\cos(x)) + (\cosh(\epsilon)-\cos(x))^{3/10}\\
    f_2(x;\epsilon) &= \frac{1}{\sqrt{\cosh(\epsilon)-\cos(x)}}\\
    f_3(x;\epsilon) &= \frac{\cos^2(6x)}{\sqrt{\cosh(\epsilon)-\cos(x)}}\\
    f_4(x;\epsilon) &= \frac{\sqrt{\cosh(1)+\cos(x)}}{\sqrt{\cosh(\epsilon)-\cos(x)}}.
\end{align}
Because $\cosh(\epsilon) = \cos(\epsilon i)$, all of these integrands have branch cuts extending vertically to infinity from $x=2\pi k \pm \epsilon i$, for $k\in\mathbb{Z}$. The first three integrands have no singularities other than these branch cuts, while $f_4$ has additional branch cuts extending to infinity from $2\pi k + \pi\pm 1i$.  The results, displayed in Fig. \ref{fig:many_periodic_examples}, suggest that the iterated sine mapping method usually reaches machine precision with a similar or smaller number of quadrature nodes $n$ compared to the Jacobi amplitude mapping, and that these two methods give much better results than the other three methods. 
We now make more detailed comments about each row of the figure. 

For the first integrand $f_1$, a sum of logarithmic and fractional-power singularities, all five methods converge with the predicted slopes (top row of Fig.  \ref{fig:many_periodic_examples}). The Jacobi amplitude mapping has the best slope, but the iterated sine mapping reaches machine precision at the same time or slightly earlier, with about $n=50$ quadrature nodes. 

For the second integrand $f_2$, which has an inverse root singularity, the Jacobi amplitude mapping gives remarkably good results. This is something of a lucky accident particular to this integrand (if one multiplies $f_2$ by $\cos(x)$, the result, not pictured, is similar to the top row of Fig. \ref{fig:many_periodic_examples}). We also see that the boundary correspondence mapping converges twice as quickly as predicted when $n$ is odd, but converges at the expected rate for even $n$ (the precise $n$ sampled in Fig. \ref{fig:many_periodic_examples} are the multiples of 7 up to 147). For the smaller values of $B$ (middle and right column of Fig. \ref{fig:many_periodic_examples}), this is still much slower than the ISM or JAM results. 

The third integrand $f_3$ is the product of $f_2$ and the entire function $\cos^2(6x)$, which grows rapidly away from the real line. Therefore we expect the constant $M$ in Theorem \ref{thm:trap} to play a more influential role in the third row of Fig. \ref{fig:many_periodic_examples} than in the second; in particular, a larger $n$ will be required before the errors decrease at the predicted rate. When $\epsilon=1/1000$, both the iterated sine mapping and the Jacobi amplitude mapping converge more slowly than predicted, and the iterated sine mapping significantly outperforms the Jacobi amplitude mapping. We also remark that the splitting method, while not converging as quickly as the best methods, does converge at the predicted rate, which is unsurprising given that it uses analyticity only in two thin ellipses which do not extend far from the real line. 

Finally we turn to the fourth integrand, which has a different domain of analyticity. While both the JAM and BCM methods have extensions to the the case of multiple singularities \cite{tee2006adaptive,berrut2020periodic}, we choose to construct the mappings with reference only to the singularities at $0\pm \epsilon i$, ignoring the more distant ones at $\pi \pm i$. 
This means that the results based on conformal mapping should converge somewhat more slowly than expected, a prediction confirmed by the last row of Fig. \ref{fig:many_periodic_examples}. However, for small $B$ the iterated sine mapping is much less impaired by this failure of analyticity than the Jacobi amplitude and boundary correspondence mappings. We can explain this difference by examining the illustrations of the three mappings in Fig. \ref{fig:cmp}. Indeed, on the lines $x=\pm\pi+\xi i$, the iterated sine mapping assumes analyticity only for a bounded range of $\xi$, while the other mappings assume analyticity for all $\xi$. This makes the singularity at $\pm \pi \pm1i$ more hazardous for the BCM and JAM methods, which rely more completely on the analyticity we have assumed. 

\subsection{Aperiodic problems with nonreal singularity} 
For this setting we modify the test problems slightly so that the integration interval is $[-1,1]$ and the singularity lies at $2/3 \pm \epsilon i$. Specifically, we use 
\begin{align}
\label{eq:firstg}
    g_1(x;\epsilon) &= -\log(\cosh(x-2/3)-\cos(\epsilon)) + (\cosh(x-2/3)-\cos(\epsilon))^{0.3}\\
    g_2(x;\epsilon) &= \frac{1}{\sqrt{\cosh(x-2/3)-\cos(\epsilon)}}\\
    g_3(x;\epsilon) &= \frac{\cos^2(6\pi x)}{\sqrt{\cosh(x-2/3)-\cos(\epsilon)}}\\
    g_4(x;\epsilon) &= \frac{\sqrt{\cosh(x+2/3)-\cos(1)}}{\sqrt{\cosh(x-2/3)-\cos(\epsilon)}}.
\label{eq:lastg}
\end{align}
Because the integration interval is smaller by a factor of $\pi$, we take $\epsilon\in \{1/30$, $ 1/300$, $1/3000\}$ to obtain problems of comparable difficulty to the previous subsection. 

\begin{figure}
    \centering
    \includegraphics[width=\linewidth]{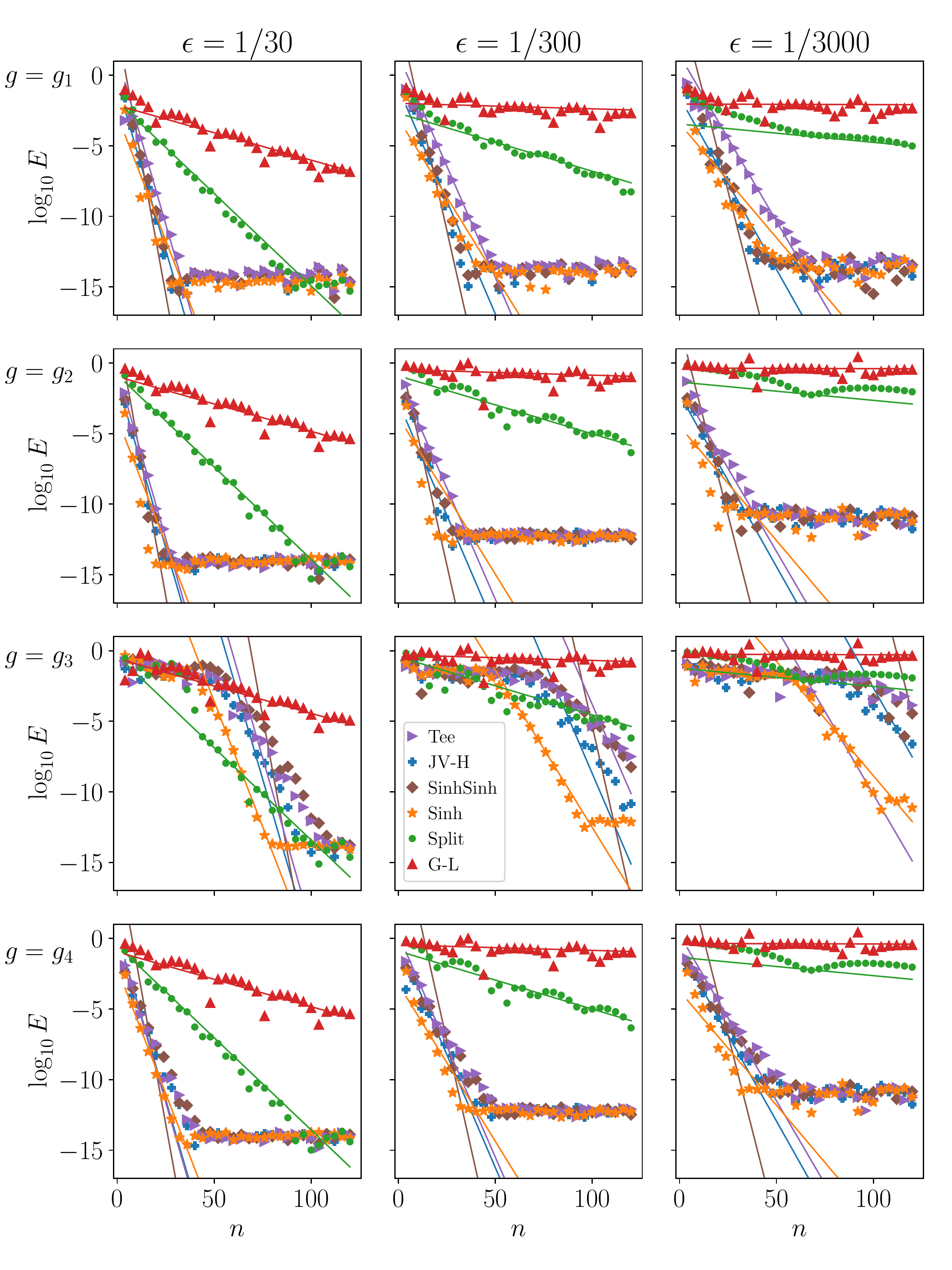}
    \caption{Integration of the functions \eqref{eq:firstg} - \eqref{eq:lastg} by the methods of Section \ref{sec:apc}. Here ``Tee'' refers to the Jacobi elliptic mapping, ``JV-H'' is the composite mapping introduced by Jafari-Varzaneh and Hosseini, ``SinhSinh'' and ``Sinh'' are the iterated and ordinary hyperbolic sine mappings, ``Split'' is the decomposition strategy, and ``G-L'' is standard Gauss-Legendre quadrature.  Among these, the ordinary sinh method generally gives the best results. }
    \label{fig:many_np_c_examples}
\end{figure}

As above, the third integrand is the product of the second integrand and an entire function, while the fourth integrand has an additional pair of singularities at $x=-2/3\pm i$.  The results, given in Fig. \ref{fig:many_np_c_examples}, are very similar to the periodic examples. In particular, the mapping that uses the Jacobi elliptic functions to carry a Bernstein ellipse onto the full domain of analyticity of the integrand does not give the best results, even though its predicted convergence rate is the best; instead, the more cautious hyperbolic sine mapping appears to be the best general choice. For these problems, the iterated sinh map does not converge at the large rate suggested by \eqref{eq:sinhsinhtstar}; its performance is similar to Tee's elliptic mapping. 

All of these examples have $A=2/3$. In other tests, not plotted, we found similar results for integrands with singularities at $1+\epsilon + \epsilon i$ for $\epsilon = 1/30, 1/300, 1/3000$.  In particular, the sinh map again gave the best results overall.

\subsection{Aperiodic problems with real singularity} 
We now integrate on $[-1,1]$ with singularity at $1 + \epsilon \in\mathbb{R}$. Specifically, we use 
\begin{align}
\label{eq:firsth}
    h_1(x;\epsilon) &= -\log(1+\epsilon-x) + (1+\epsilon-x)^{0.3}\\
    h_2(x;\epsilon) &= \frac{1}{\sqrt{1+\epsilon-x}}\\
    h_3(x;\epsilon) &= \frac{\cos^2(6\pi x)}{\sqrt{1+\epsilon-x}}\\
    h_4(x;\epsilon) &= \frac{\sqrt{\cosh(x+2/3)-\cos(1)}}{\sqrt{1+\epsilon-x}}.
\label{eq:lasth}
\end{align}
We again take $\epsilon\in \{1/30, 1/300, 1/3000\}$.

\begin{figure}
    \centering
    \includegraphics[width=\linewidth]{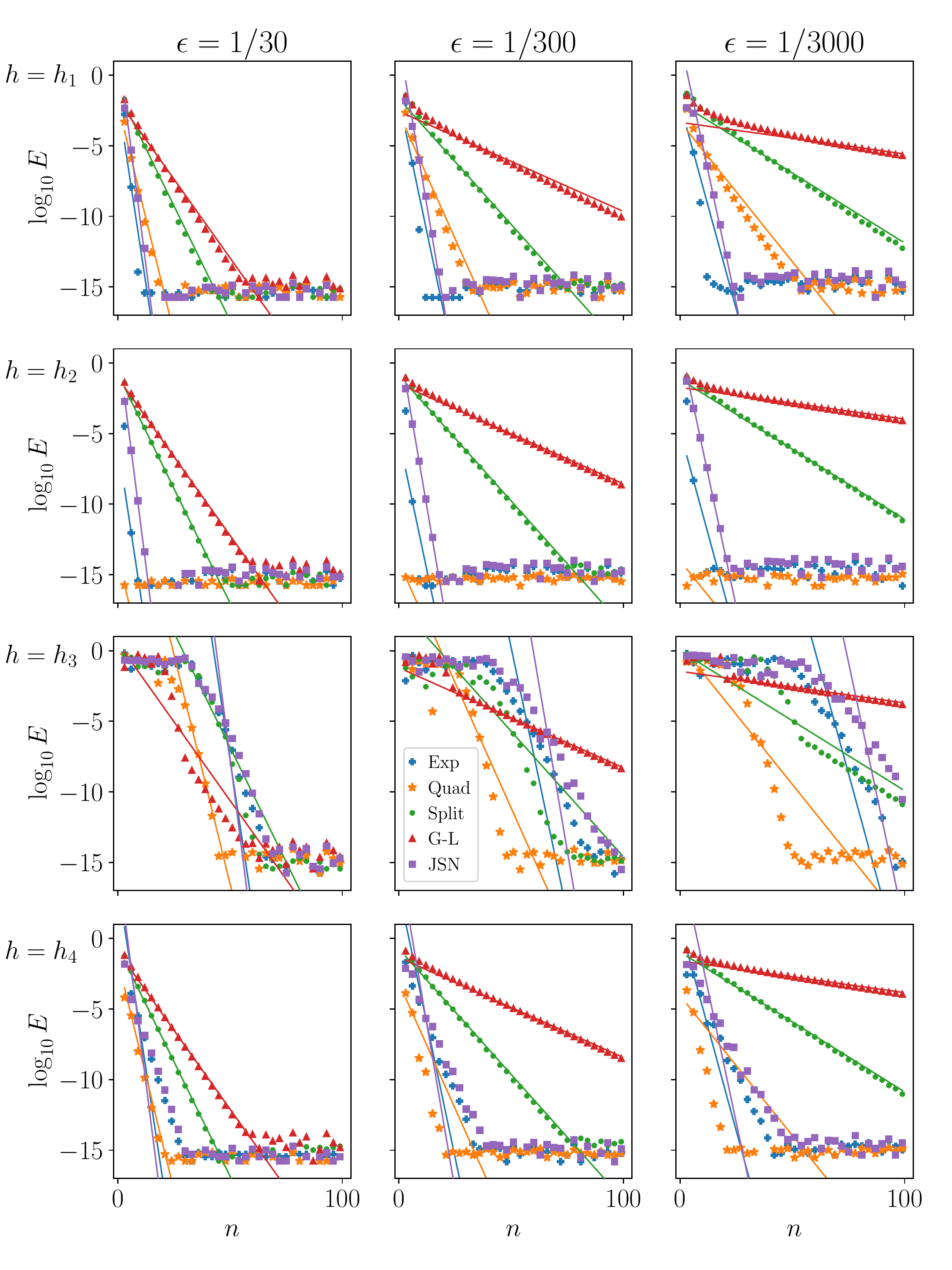}
    \caption{Integration of the functions \eqref{eq:firsth} - \eqref{eq:lasth} the exponential, quadratic and Jacobi elliptic sine maps along with the decomposition method and ordinary Gauss-Legendre quadrature as discussed in Sec. \ref{sec:apr}. The quadratic transformation is best in the last three rows, while the exponential map is best in the first row. }
    \label{fig:many_np_r_examples}
\end{figure}

The results, in Fig. \ref{fig:many_np_r_examples}, indicate that the quadratic transformation gives the best results for $h_2$, $h_3$, and $h_4$, while the exponential transformation is better for $h_1$. 

\section{Application: evaluation of single-layer potentials in Stokes flow} \label{sec:SLP}

\begin{figure}
    \centering
\includegraphics[width=\linewidth]{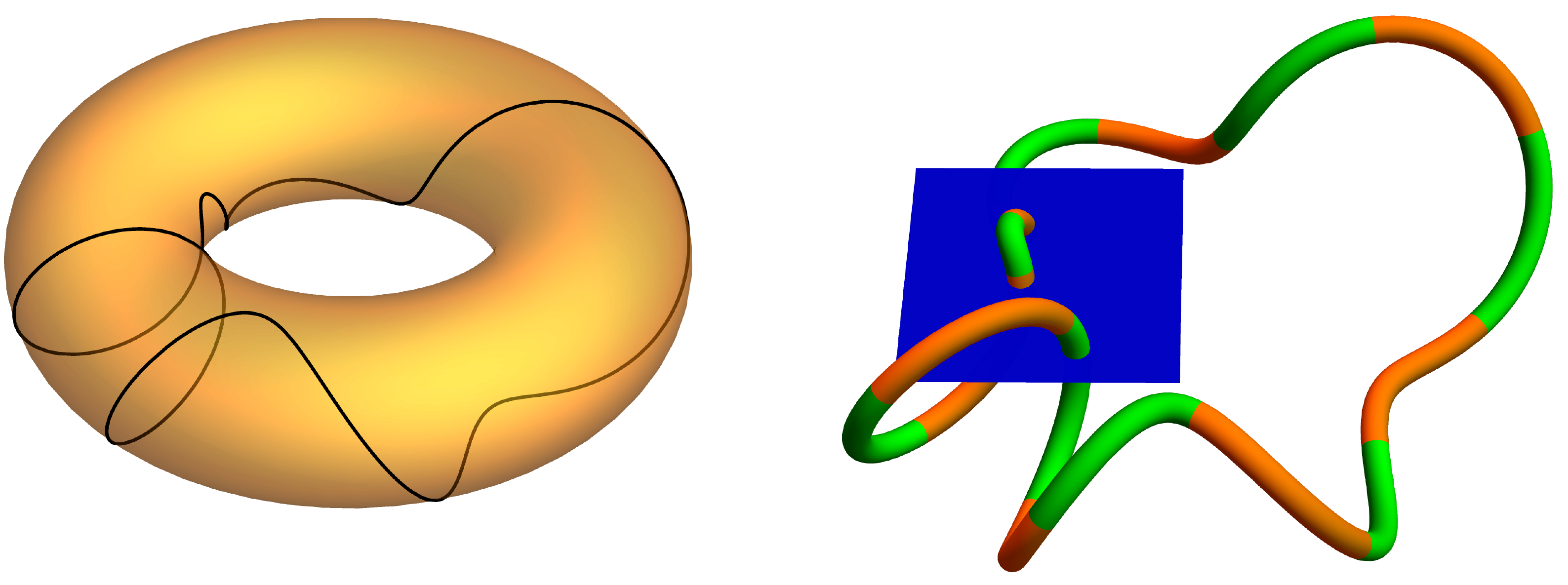}
    \caption{To demonstrate our integration strategies we consider the problem of evaluating a Stokes single-layer potential where the target point is near but not on the surface of integration. The surface is a fiber whose centerline follows a closed path on the surface of a torus (left). We then let the target point range over a square domain that is punctured in three places by the fiber (right). The integrals become numerically challenging when the target point approaches the surface. To organize the integration we subdivide the fiber surface into sixteen panels, as depicted in green and orange at right. }
    \label{fig:fiber}
\end{figure}

As an application of the preceding methods, we will evaluate some nearly singular surface integrals that arise in the study of viscous fluid flow. Specifically, we will consider the Stokes single-layer potential defined by 
\begin{equation}
\mathcal{S}(\bm x) = \int_D \left(\frac{\bm f(\bm y)}{\lvert \bm x-\bm y \rvert } + \frac{\bm x - \bm y}{\lvert \bm x-\bm y \rvert ^3}((\bm x - \bm y)\cdot \bm f(\bm y))\right) \, dS_{\bm y}
\label{eq:slp}
\end{equation}
where $D$ is the surface of the slender fiber depicted in the right panel of Fig. \ref{fig:fiber}. The integral has a physical meaning: it is the velocity field that results when point forces of strength $\bm f$ are distributed over the surface $D$. 
The integral is challenging when the \emph{observation point} $\bm x$ is near but not on the surface $D$. Although the surface has an irregular shape, the integral has a known exact solution when the density $\bm f$ is equal to the outward pointing surface normal: in this case the resulting velocity is zero everywhere, and the surface traction $\bm f$ is due solely to hydrostatic pressure. This arrangement allows us to test various quadrature strategies in a setting where the geometry is nontrivial but an exact solution (zero) is known. Note that we have omitted the usual factor of $1/(8\pi)$ in \eqref{eq:slp}. 

The problem of computing nearly singular surface integrals arises in many situations when a linear PDE is solved via integral equations. For the Stokes PDE, a number of methods have been developed that take into account the particular form of the Stokes fundamental solutions. Prominent recent examples include singularity subtraction  and local expansion \cite{tlupova2013nearly,tlupova2019regularized}, quadrature by expansion \cite{af2016fast}, and singularity swapping \cite{af2021accurate}. Here we attempt to address this challenge without using special knowledge about the particular form of the integrands. 

To describe the fiber surface, we begin by parameterizing the surface of a torus whose centerline has unit radius and whose circular cross sections have radius $0.4$: 
\begin{equation}
    \bm v(\theta,\phi) = (1+0.4\cos(\phi))\begin{pmatrix}
    \cos\theta\\\sin\theta\\0\end{pmatrix}+0.4\sin\phi\begin{pmatrix}
    0\\0\\1\end{pmatrix}.
\end{equation}
Then we construct a closed curve on this surface by letting $\theta$ and $\phi$ be periodic functions of a (non-arclength) parameter $s$: 
\begin{equation}
    \bm w(s) = \bm v\Big(s,2 \exp(\cos(s + 1)) \cos(2 s) + 2 s\Big).
\end{equation}
The path $\bm w$ appears on the surface of the torus in the left panel of Fig. \ref{fig:fiber}. 
Let $\bm T(s)$ denote the unit tangent vector for the curve, $\bm T(s) = \bm w'(s) / \lvert\bm w'(s)\rvert$. To define the fiber surface we need vectors $\bm N(s)$ and $\bm B(s)$ so that $\{\bm T, \bm N,\bm B\}$ is an orthonormal frame. To avoid the derivatives involved in the standard Frenet definition, we plotted $\bm T(s)$ as a path on the unit sphere and noticed that it always remains far from the line through $\pm\bm p$, where $\bm p = \langle 10,3,6\rangle$. Thus we can complete the frame by putting 
\begin{equation}
    \bm N(s) = \frac{\bm T(s) \times \bm p}{ \lvert\bm T(s) \times \bm p\rvert},\qquad
\bm B(s) = \bm T(s) \times \bm N(s).
\end{equation}
Finally we define the fiber surface $D$ via 
\begin{equation}
    \bm y(s,t) = \bm w(s) + \epsilon\cos(t)\bm N(s) + \epsilon\sin(t)\bm B(s)
\end{equation}
where $s$ and $t$ both range from $0$ to $2\pi$, and the radius $\epsilon$ is a constant (we take $\epsilon = 0.05$).
The surface normal vector for the fiber is 
\begin{equation}
    \bm f(s,t) = \bm \nu(s,t) = \cos(t)\bm N(s) + \sin(t)\bm B(s)
\end{equation}
while the Jacobian or surface integration weight is 
\begin{equation}
    J(s,t) = \left\lvert\frac{\partial \bm y}{\partial s} \times \frac{\partial \bm y}{\partial t}\right\rvert = 0.1 (\lvert\bm w'(s)\rvert-0.1 \kappa_1 \cos(t) -0.1\kappa_2\sin(t))
\end{equation}
where $\kappa_{1,2}$ are defined by $\kappa_1 = \bm w'(s)\cdot \bm N(s)$ and $\kappa_2 = \bm w'(s)\cdot \bm B(s)$. 

\subsection{Reference solution}
We first evaluate the integral using ordinary Gauss-Legendre and trapezoidal quadrature. 
We subdivide the fiber surface into 16 panels using a heuristic that considers the panel lengths and their maximum curvatures; the panel endpoints are at $s\in\{ 0$, $0.58$, $ 1.21$, $ 1.83$, $2.35$, $2.76$, $3.19$, $3.86$, $4.26$, $4.65$, $5.04$, $5.24$, $5.41$, $5.57$, $5.75$, $ 6.05$, $2\pi\}$. 
On each panel, we use the outer product of a Gauss-Legendre grid in $s$ and an equally spaced grid in $t$, using the same quadrature rule for every target point. The target points range over a square domain $\{(x,y,0):0<x<1, -5/8<y<3/8\}$ that is punctured in three locations by the fiber surface. The exact velocity is zero, and we depict the norm of the computed velocity at each target point as a contour plot in the first column of Fig. \ref{fig:slpcontours}.  Predictably, we see that this combination of Gauss-Legendre and trapezoidal rules is effective only when the target is far from the fiber surface. 
\subsection{Finding the singularities}
In order to improve on the reference solution using the methods described in this paper, we need to find the singularities of the inner (periodic) integrand in the complex $t$-plane for fixed $s$, as well as the singularities of the outer integrand in the complex $s$-plane. Although the outer integral in $s$ is also $2\pi$-periodic, we chose to integrate separately on each of the panels, leading to sixteen aperiodic subproblems. 

We can find the singularities for the periodic, inner problem on paper as follows. 
We begin by writing \eqref{eq:slp} as a double integral. With $\bm r(s,t) = \bm x - \bm y(s,t)$ we have 
\begin{equation}
\mathcal{S}(\bm x) = \int_0^{2\pi}\int_0^{2\pi} \left(\frac{\bm f(\bm y(s,t))}{\lvert\bm r(s,t)\rvert} + \frac{\bm r(s,t)}{\lvert\bm r(s,t)\rvert^3}(\bm r(s,t)\cdot \bm f(\bm y(s,t)))\right) J(s,t)\, dt\,ds.
\label{eq:slpit}
\end{equation}
For fixed $s$, the vector $\bm r(s,t)$ traces out a circle at constant speed as $t$ varies. The denominators can therefore be written as trigonometric functions of $t$ and we can solve analytically for the complex value of $t$ where they vanish. 
To do this we write 
\begin{align*}
  \lvert\bm r(s,t)\rvert^2 &= \big\lvert\bm x - \bm w(s) - \epsilon\cos(t)\bm N(s) - \epsilon \sin(t)\bm B(s)\big\rvert^2\\
    &=\lvert\bm x - \bm w(s)\rvert^2 + \epsilon^2 -2\epsilon\Big(\cos (t) (\bm x-\bm w(s))\cdot\bm N(s) + \sin(t)(\bm x-\bm w(s))\cdot\bm B(s)\Big)\\
    &= \lvert\bm x - \bm w(s)\rvert^2 + \epsilon^2 -{2\epsilon}{\sqrt{((\bm x-\bm w(s))\cdot\bm N(s))^2+((\bm x-\bm w(s))\cdot\bm B(s))^2}}\cos(t-\xi) 
\end{align*}
where the angle $\xi$ is defined by 
\begin{align*}
\cos(\xi) &= \frac{(\bm x-\bm w(s))\cdot\bm N(s)}{\sqrt{((\bm x-\bm w(s))\cdot\bm N(s))^2+((\bm x-\bm w(s))\cdot\bm B(s))^2}},\;\\
\sin(\xi) &= \frac{(\bm x-\bm w(s))\cdot\bm B(s)}{\sqrt{((\bm x-\bm w(s))\cdot\bm N(s))^2+((\bm x-\bm w(s))\cdot\bm B(s))^2}}.
\end{align*}
Therefore, $\bm r(s,t)$ vanishes when 
\[
\cos(t-\xi) = \frac{\lvert\bm x - \bm w(s)\rvert^2 + \epsilon^2}{2\epsilon\sqrt{((\bm x-\bm w(s))\cdot\bm N(s))^2+((\bm x-\bm w(s))\cdot\bm B(s))^2}}
\]
and we find that the required value of $t$ is
\begin{equation}
t^* = \xi + \sqrt{-1}\arccosh\left(\frac{\lvert\bm x - \bm w(s)\rvert^2 + \epsilon^2}{2\epsilon\sqrt{((\bm x-\bm w(s))\cdot\bm N(s))^2+((\bm x-\bm w(s))\cdot\bm B(s))^2}}\right).
\label{eq:innertstar}
\end{equation}
This allows us to choose a quadrature for the inner integral using knowledge of the integrand's complex singularity, as in Section 2. 

We will need some numerical rootfinding to locate the singularities of the outer integral. We note that the inner integral will diverge if the imaginary part of $t^*$ in \eqref{eq:innertstar} vanishes. Therefore, the complex singularities of the outer integrand are the solutions of the equation
\begin{align}
    \Big(\lvert\bm x - \bm w(s)\rvert^2 + 0.1^2\Big)^2 = 0.04\Big(((\bm x-\bm w(s))\cdot\bm N(s))^2 +((\bm x-\bm w(s))\cdot\bm B(s))^2\Big). \label{eq:outersing}
\end{align}
To obtain these roots we find the eigenvalues of a $51\times 51$ Chebyshev colleague matrix. For each computed root we then find the corresponding Bernstein ellipse parameter $\rho$. If all computed $\rho$-values are greater than 2, we revert to ordinary Gauss-Legendre quadrature; otherwise we use the root with the smallest $\rho$-value to accelerate the quadrature in the outer integral.

\subsection{Improvement via decomposition and conformal mapping}
\begin{figure}
    \centering
    \includegraphics[width=\linewidth]{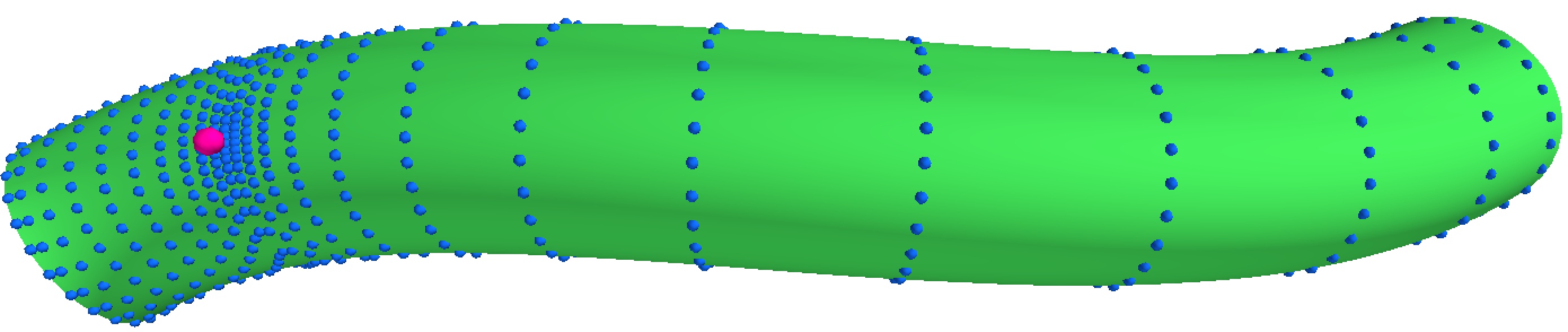}
    \caption{A surface quadrature rule with $24\cdot24$ nodes (blue dots) on one panel, customized for a target point (pink dot) that is near but not on the surface. We generated this rule with the hyperbolic sine transformation in the centerline direction, $s$, and the iterated sine map in the circumferential direction, $t$ (calculated separately for each value of $s$). }
    \label{fig:dotsontube}
\end{figure}

Overall, the surface quadrature procedure that we have outlined is laborious: for each target point and within each panel, we find a customized 1D quadrature rule for the aperiodic outer integral. Then, for each of the resulting outer quadrature nodes, we find a customized 1D quadrature rule for the periodic inner integral. An example of the resulting surface quadrature rule appears in Fig. \ref{fig:dotsontube}. Although our goal is to demonstrate the accuracy of the underlying quadratures rather than to address the fast generation of rules for many target points, we make one adjustment for the sake of efficiency: for each panel, we revert to the reference Gauss-Legendre / trapezoid scheme for all targets whose distance to the panel surface is greater than $7\epsilon$. This allows us to use the same quadrature rule simultaneously for many distant targets. 

\begin{figure}
    \centering
    \includegraphics[width=\linewidth]{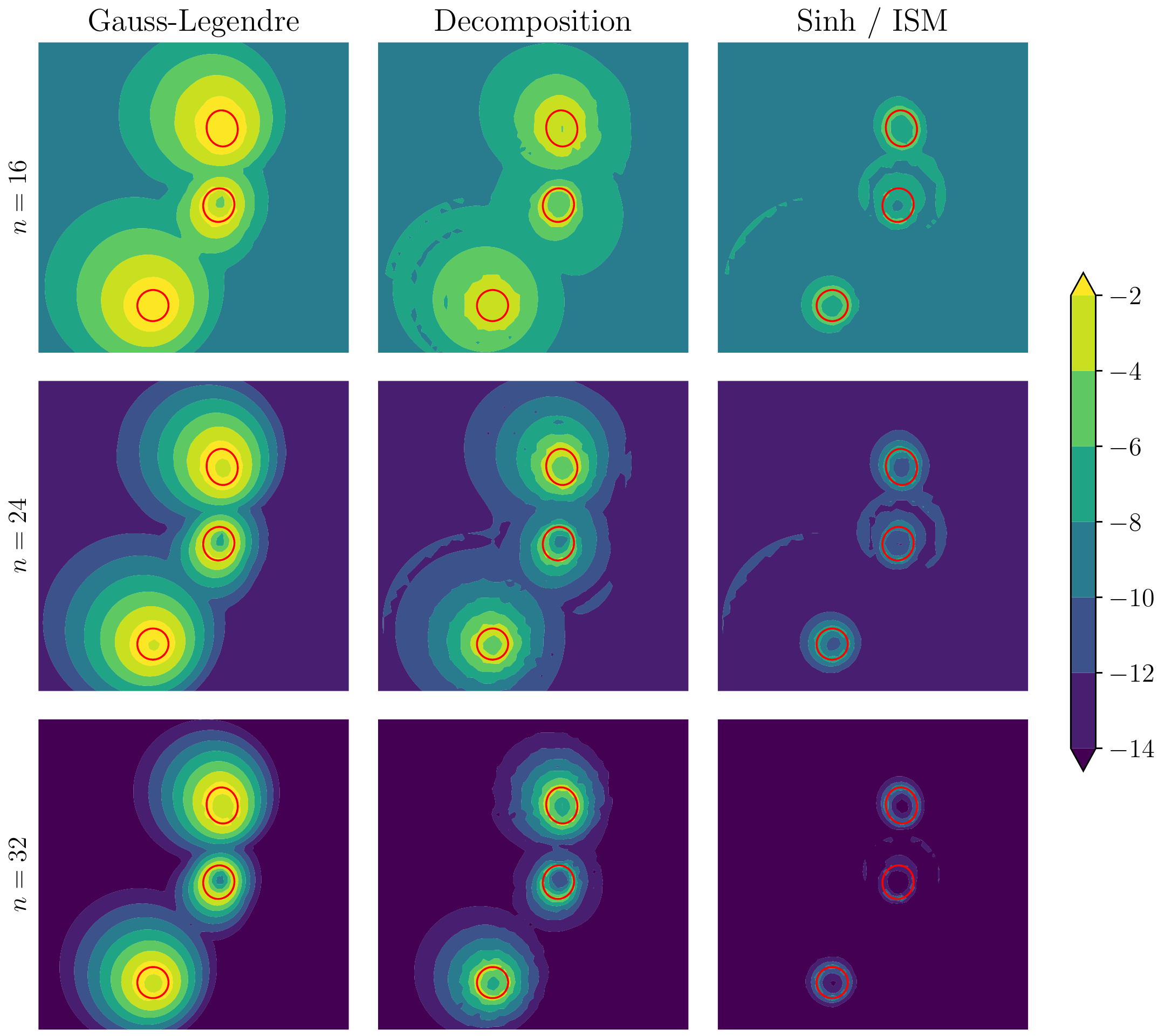}
    \caption{Contour plot of errors in the single-layer problem as the target point ranges over a square domain, punctured in three places by the slender fiber. The error bar reports $\log_{10}(\lvert\bm u(\bm x)\rvert)$ where $\bm u(\bm x)$ is the single-layer velocity induced at the target point $\bm x$ by a surface traction which equals the unit vector ($\bm u=\bm0$ if the integration is done correctly). The columns of the figure show different quadrature strategies, while the rows of the table show differing densities of quadrature nodes; for example, $n=32$ means that we use $32\cdot32$ nodes on each panel (there are always 16 panels).   
    }
    \label{fig:slpcontours}
\end{figure}

The second column of Fig. \ref{fig:slpcontours} shows the result of combining the periodic decomposition method for the inner integral with the aperiodic decomposition method in the outer integral. This leads to more correct digits than the reference solution when the target point is near the surface. 
However, the third column, the result of using the iterated sine map for the inner integral and the hyperbolic sine map for the outer integral, is dramatically better. We chose the sinh and ISM methods  because of their simplicity and superior performance on the tests of section \ref{sec:examples}. However, other combinations of the conformal mapping strategies, including the Jafari-Varzaneh map and the various options based on Jacobi elliptic functions, yield similar (but not better) results for this application. Note that we did not use the results of \ref{sec:apr} because the numerical rootfinding procedure always gave a nonzero imaginary part. 

Our Fig. \ref{fig:slpcontours} should be compared with Figures 2, 3, and 8 of \cite{af2021accurate}, a related work where the authors make use of additional information about the nature of the singularities, not merely their location, and accordingly develop a more powerful but less general method for accelerating the quadrature of the nearly singular integrals. We note that they developed methods for finding and using multiple root pairs, but found in practice that this did not give better results than the methods based on a single root pair. 

\section{Conclusion}
\label{sec:conclusion}
We surveyed a number of possible strategies for accelerating the quadrature of nearly single integrals given knowledge of the location of the nearby singularity. We found that the splitting methods are less effective than the conformal maps. Among the many possible conformal maps, we can make suggestions for general use: we recommend the iterated sine map for periodic problems, the sinh map for aperiodic problems with complex singularity, and the quadratic map for aperiodic problems with real singularity.

A natural priority for future work is to ask what additional improvements are possible when we have information about the nature of the nearest singularity in addition to its location. The integrand $f(x) = ((x-0.3)^2+\epsilon^2)^p$ has the same domain of analyticity for $p=2.5$ and $p=-2.5$, and accordingly the Gauss-Legendre quadrature errors will \emph{eventually} decay at the same rate in either case. However, it is also true that Gauss-Legendre integration reaches machine precision much more quickly for the version with $p>0$. This statement that can be made more precise with the aid of the theorems of Sec. \ref{sec:thms}, but it remains to make use of this information to optimize the choice of conformal map.

A second issue arises when many integrals over the same surface, but with varying singularities, need to be computed. The strategy described here would use customized quadratures for each of the singularities, resulting in a possibly expensive interpolation operation. It would be valuable to partition the complex plane into regions such that any singularity within a region can be integrated to prescribed accuracy with a precomputed reference set of quadrature nodes and weights.

\subsection{Acknowledgements}
We thank Nick Trefethen for helpful feedback on a draft. We acknowledge support from the National Science Foundation (DMS-1907796) and from Macalester College.

\bibliographystyle{spmpsci}      
\bibliography{references}

\end{document}